\documentclass[onecolumn,a4paper,draftcls]{article}
\usepackage{dsfont}
\usepackage{amsmath, amssymb, amsthm}
\usepackage{adjustbox}
\usepackage{graphicx}
\usepackage{enumerate}
\usepackage{lineno}
\usepackage{hyperref}
\usepackage{xcolor}
\usepackage{setspace}
\def\Blue#1{\textcolor{blue}{#1}}

\newtheorem{remark}{Remark}[section]

\newtheorem{definition}{Definition}[section]

\newtheorem{assumption}{Assumption}[section]

\newtheorem{problem}{Problem}[section]

\newtheorem{theorem}{Theorem}[section]

\newtheorem{corollary}{Corollary}[section]

\newtheorem{lemma}{Lemma}[section]

\newtheorem{proposition}{Proposition}[section]

\def\R{\mathbb{R}}
\def\N{\mathbb{N}}
\def\K{\mathcal{K}}
\def\L{\mathcal{L}}

\def\cal{\mathcal}

\title{\LARGE \bf
On the robustness of self-triggered sampling of nonlinear control systems.}

\author{U. Tiberi, K.H. Johansson, Fellow,~IEEE
\thanks{K.H. Johansson is with ACCESS Linneaus Center, KTH Royal Institute of Technology, Stockholm, Sweden. Email: kallej@kth.se}
\thanks{U. Tiberi is with Volvo Group Trucks Technology, G\"oteborg, Sweden. Email: ubaldo.tiberi@volvo.com}
}
\providecommand{\keywords}[1]{\textbf{\textit{Index terms---}} #1}

\begin{document}
 \maketitle                                

\begin{abstract}
We address robustness issues of self-triggered sampling with respect to model uncertainties, and propose a robust self-triggered sampling method. The approach is compared with existing methods in terms of sampling conservativeness and closed-loop system performance. The proposed method aims at fulfilling the gap between the event and the self-triggered sampling paradigms for what concerns robustness with respect to model uncertainties, and it generalizes most of the existing self-triggered samplers implemented up to now.
\end{abstract}


%

\section{Introduction}\label{sec:introduction}

To cope with common drawbacks raised by periodic sampling in modern control systems, such as network utilization in networked control systems~\cite{HESPANHA-IEEE-2007} or processor utilization in multi-task programming~\cite{BUTTAZZO-BOOK}, two novel sampling methods, referred to as event-based and self-triggered sampling, has been recently introduced,~\cite{ASTROM-CDC-2002}--\nocite{lunze10,GARCIA-TAC-2014,KALLE-CDC-2009,TABUADA-TAC-2007,wang11,MAZO-AUTOMATICA-2010,ANTA-TAC-2010}\cite{TIBERI-NAHS-2013}. Roughly speaking, event-based sampling consists in monitoring the system output for all the time and to update the control signal only when some event is detected, whereas self-triggered sampling consists in \emph{predicting} the event occurrence based on a system model and on the current system output. 
It has been showed that both approaches usually leads to an efficient utilization of shared resources without deteriorating the closed-loop performance.
Nevertheless, they exhibit profound differences.
Event-based methods take decisions upon the detection of an event and they can be thus categorized as \emph{reactive} methods; on the contrary, self-triggered methods are \emph{proactive} as they provide the next event occurrence time in advance. 
A notable benefit in event-based methods is that they seldom requires a model of the plant, but the event occurrences are often determined only from the output measurements, whereas in self-triggered methods an accurate system model is generally required. 
Clearly, if the model is not sufficiently accurate, the closed-loop performance under self-triggered sampling may deteriorate or, in some cases, the closed-loop system may even become unstable.

To the best of our knowledge, the problem of self-triggered control robustness versus parameter uncertainty for nonlinear systems has been little investigated, and existing methods exhibit severe limitations~\cite{DIGENNARO-IJC-2013},\cite{TIBERI-PHDTHESIS}. For instance, the approach proposed in~\cite{DIGENNARO-IJC-2013} relies on assumptions that hold only for a very narrow class of systems, thus limiting the applicable cases. Nevertheless, if such assumptions are relaxed, then both the cited methods can only guarantee a safety property of the closed-loop system which is weaker than common stability properties such as asymptotic stability of ultimate boundedness.

In contrast, our method requires milder assumptions compared to the cited work, which extends the applicability to a larger number of cases. Inspired by the Lebesgue sampling rule~\cite{ASTROM-CDC-2002}, our approach ensures uniform ultimate boundedness or, in some cases even asymptotic stability. In this note, we address both the local and the global stability cases. Finally, the proposed approach is compared with existing methods in terms of conservativeness of the sampling intervals and closed-loop performance.


%
\section{Notation and preliminaries}\label{sec:notation}
The set of natural numbers is denoted with $\N$.
The set of real numbers is denoted with $\R$, the set of positive real numbers with $\R^+$ and the set of nonnegative real numbers with $\R^+_0$, i.e. $\R_0^+=\R^+\cup \{0\}$. 
The notation $\|v\|$ is used to indicate the Euclidean norm of a vector $v\in \mathbb R^{n}$ and $\mathcal B_{r}$ indicates the closed ball centered at the origin and radius $r$, i.e. $\mathcal B_{r}=\{v:\|v\|\le r \}$. Given a set $\mathcal D$, we denote its power set with $2^{\mathcal D}$.
 Given a signal $s: \R^{+}  \to \mathbb R^{n}$, $s_{k}$ denotes its realization at time $t=t_{k}$, i.e. $s_{k}:=s(t_{k})$.
 A function $h:\cal D_p \times \cal D_q  \to \R^n$ is said to be \emph{Lipschitz continuous over $\cal D_p \times \cal D_q$} if $\|h(p_1,q)-h(p_2,q)\| \le L_{h,p}\|p_1-p_2\|$ for some $L_{h,p}>0$ and for all $p_1,p_2 \in \cal D_p, q \in \cal D_p$ and $\|h(p,q_1)-h(p,q_2)\| \le L_{h,q}\|q_1-q_2\|$ for some $L_{h,q}>0$ and for all $q_1,q_2 \in \cal D_q, p \in \cal D_p$. The constants $L_{h,p}$ and $L_{h,q}$ are called \emph{Lipschitz constant of $h$ with respect to $p$} and \emph{Lipschitz constant of $h$ with respect to $q$}, respectively.
A continuous function $\alpha:[0,a) \to +\infty,a>0$ is said to belong to class $\cal K$ if it is strictly increasing and $\alpha(0)=0$. If, in addition, $a=\infty$ and $\alpha(r) \to +\infty$ for $r \to +\infty$, then $\alpha$ is said to be of class $\cal K_{\infty}$.
A continuous function $\beta:[0,a)\times [0,\infty)\to [0,\infty)$ is said to belong to class $\K\L$ is, for each fixed $s$, the mapping $\beta(r,s)$ belongs to class $\K$ with respect to $r$ and, for each fixed $r$, the mapping $\beta(r,s)$ is decreasing with respect to $s$ and $\beta(r,s) \to 0$ as $s\to \infty$.
Given a system $\dot \xi=f(t,\xi)$, $\xi \in \mathbb R^{n}$, $\xi(t_{0})=\xi_{0}$, $f:\mathbb{R}^{+} \times \mathcal{D} \rightarrow \mathbb{R}^{n}$, where $f$ is Lipschitz continuous with respect to $x$ and piecewise continuous with respect to $t$, and where $\mathcal{D} \subset \mathbb{R}^{n}$ is a domain that contains the origin, we say that the solutions are UUB if there exists three constants $a,b,T>0$ independent of $t_{0}$ such that for all  $\|\xi_{0}\|\le a$ it holds $\|\xi(t)\| \le b$ for all $t \ge t_{0}+T$, and globally UUB (GUUB) if $\|\xi(t)\| \le b$ for all $t \ge t_{0}+T$ and for arbitrarily large $a$. The value of $b$ is referred as the \emph{ultimate bound}.
%
%
%
%
\section{System architecture}
\begin{figure}[tbp]
  \centering
   \includegraphics[scale=0.6]{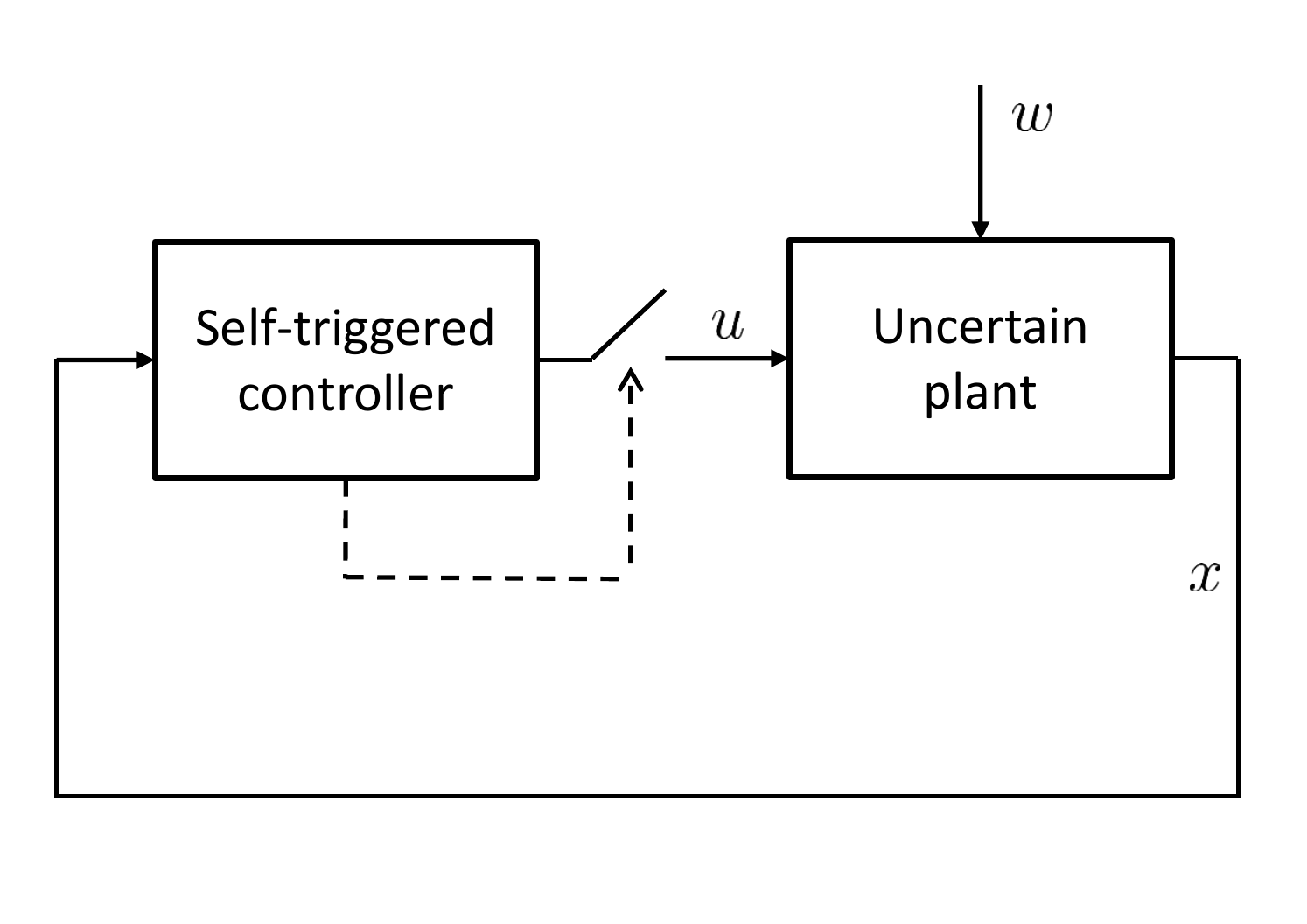}
   \caption{The control system architecture.}
   \label{fig:SysArchitecture}
\end{figure}

We consider the system architecture depicted in Fig.~\ref{fig:SysArchitecture}. 
The control system includes an uncertain plant subject to external disturbances $w$ and a self-triggered controller, i.e. a controller which computes both the new control signal and its update instant. 
The plant's dynamics are of the form
\begin{equation}\label{eqn:nominal plant}
\dot \xi = f(\eta,\xi,u,w)\,,
\end{equation}
where $f$ is Lipschitz continuous and where $\xi \in \mathcal D_{\xi}  \subseteq \R^{n_{\xi}}$ is the state vector, $u\in \mathcal D_{u} \subseteq \R^{n_{u}}$ is the input vector, $\eta $ is a vector of (possible time-varying) uncertain parameters in a compact set $\mathcal D_\eta \subset \R^{n_{\eta}}$ and $w \in \mathcal D_{w}\subseteq \R^{n_{w}}$ is a piecewise bounded external disturbance vector with bound $\|w\|\le \bar w$ 
We assume that there exists a Lipschitz continuous state feedback control law $\kappa:\mathcal D_{\xi}\to \mathcal D_{u}$ such that the closed-loop dynamics satisfying
\begin{equation}\label{eqn:closed-loop CT}
\dot \xi = f(\eta,\xi,\kappa(\xi),w)\,.
\end{equation}
are asymptotically stable for $w=0$ and UUB for all $w\in \mathcal{D}_w\backslash \{0\}$.
%
%
Our goal is to determine a function $\Gamma:\R^{2n}\to \R$ and to predict, at each time $t=t_k$\Blue{,} the time instant $t_{k+1}$ defined as
\begin{equation}\label{eqn:STC def}
t_{k+1}=t_k+\min\{t \Blue{>} t_k:\Gamma(x,x_k)=0\}\,,
\end{equation}
such that the sampled-data system
\begin{equation}\label{eqn:sampled-data}
\dot x = f(\eta,x,\kappa(x_{k}),w)\,, \quad t\in [t_{k},t_{k+1})\,,
\end{equation}
where $x \in \mathcal{D}_{\xi}$, is UUB for all $w\in \mathcal{D}_w$ and such that $t_{k+1}-t_k\ge h_{\min}$ for some $h_{\min}>0$ and for all $k$. 
Although existing self-triggered samplers may still apply for stabilizing uncertain systems of the form~\eqref{eqn:nominal plant}, the performance of the closed-loop system may not be acceptable or it may becomes unstable, as we will discuss in the next Section.
%
%
%
%
%
%
\section{Motivating example}\label{sec:motivating example}
\begin{figure}[tp]
\centering
\includegraphics[scale=0.63]{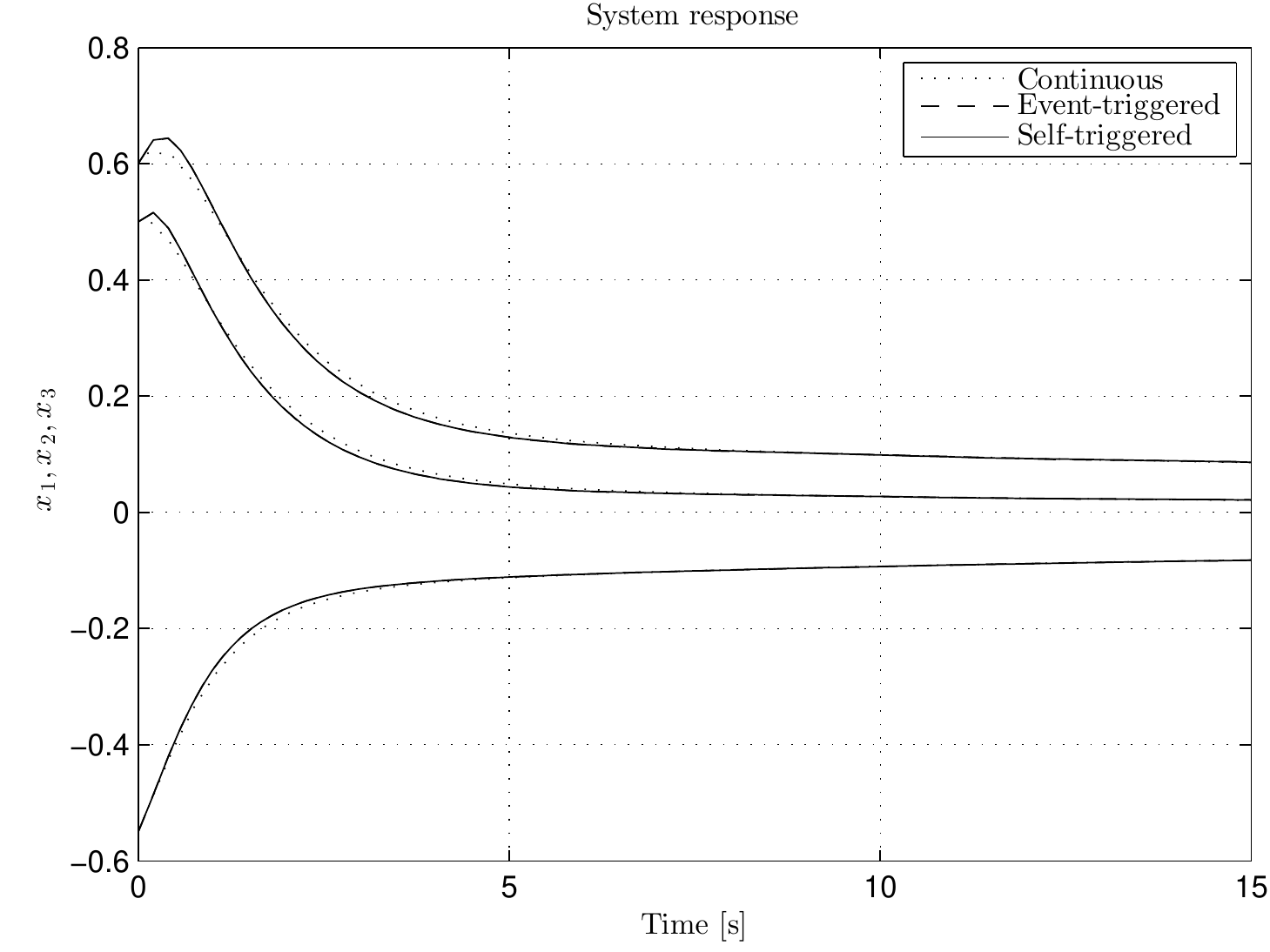}
\caption{System response with $\eta=\eta_n, \eta_n=1$}.
\label{fig:EvoEx}
\end{figure}

Consider the rigid-body control example in~\cite{ANTA-TAC-2012}, which dynamics satisfies
\begin{align}\label{eqn:rigid-body}
\dot \xi_{1}&=  u_{1}\,,\nonumber\\
\dot \xi_{2}&=  u_{2}\,,\\
\dot \xi_{3}&= \eta \xi_{1}\xi_{2}\nonumber\,,
\end{align}
and let $\eta$ to be an uncertain parameter. Let $\eta_n=1$ be the assumed value of the uncertainty for designing both the continuous controller and the self-triggered sampler.
With this setting, a control law to globally stabilize system~\eqref{eqn:rigid-body} if $\eta=\eta_n$ is given by
 $u_{1}=-\xi_{1}\xi_{2}-2\xi_{2}\xi_{3}-\xi_{1}-\xi_{3}$ and $u_{2}=2\xi_{1}\xi_{2}\xi_{3}-3\xi_{3}^{2}-\xi_{2}$
and a self-triggered sampler implementation considers the sampling rule $\Gamma(x,x_k):=\|x_k-x\|^{2}-0.79^2\sigma^2\|x\|^{2}$ where  $0<\sigma<1$, see~\cite{ANTA-TAC-2012}. If for the real system it holds $\eta=\eta_n$, then the response of continuous-time, the event and the self-triggered implementation of the controller would be fairly similar as shown in Figure~\ref{fig:EvoEx}. 

\begin{figure}[tp]
\centering
\includegraphics[scale=0.63]{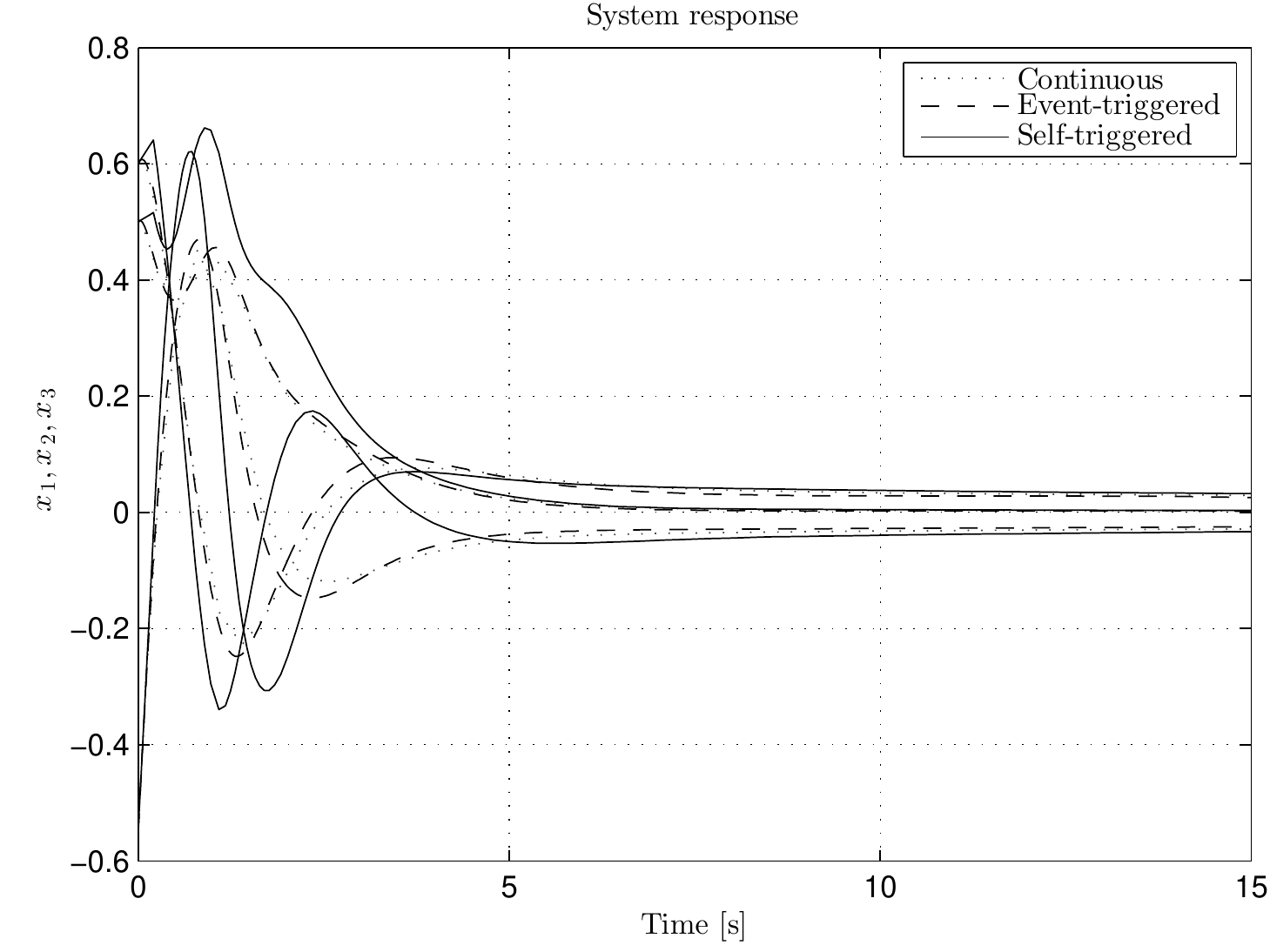}
\label{fig:TrigEx}
\caption{System response with $\eta=\eta_r, \eta_r=8$}.
\label{fig:EvoDelta}
\end{figure}

Nevertheless, assume now that for the real system dynamics it holds instead $\eta=\eta_r$, with $\eta_r=8$, and assume to use the same controller and self-triggered sampler as for the case $\eta=\eta_n$.   
As shown in Figure~\ref{fig:EvoDelta}, the closed-loop system performance deteriorates under self-triggered sampling, although both in the event and the continuous case we experience a satisfactory response. This is because in the event-triggered scheme the condition $\Gamma(x,x_k)=0$ is constantly evaluated based on a constant monitoring of the state $x(t)$, and then the time $t_{k+1}$ defined in~\eqref{eqn:STC def} are correctly determined. In the self-triggered implementation, the times $t_{k+1}$ may mismatch with the ones defined in~\eqref{eqn:STC def} since the prediction for which $\Gamma(x,x_k)=0$ is based on an imperfect model. Note that although the continuous-time controller exhibits a certain degree of robustness, this is unfortunately not enough to ensure good performance of its self-triggered implementation, but the self-triggered strategy shall also be robust.

We wish further to highlight that the event-based sampling rule implicitly defined by the function $\Gamma(x,x_k)$ in this example only represents a sufficient condition for the closed-loop stability. This means that the self-triggered implementation based on an imperfect model may not fulfill such a condition for all the time and then the closed-loop system stability is also jeopardized.

%
 %
%
%
Unfortunately, the inclusion of parameter uncertainties in the framework proposed in~\cite{ANTA-TAC-2012} does not appear to be straightforward, and leaves room for future research. Nevertheless, in this note we follow a different approach by proposing a method which applies to every robustly stabilizable nonlinear system and which ensure UUB of the sampled-data system trajectories.
 
%
%
%
%
%
%
%
%
%
\section{Robust self-triggered sampling}\label{sec:robust self-triggered sampling}
In this section we present the main result of this note. We first consider the local stabilizability and then the global stabilizability case. 
\subsection{Local analysis: exponentially stabilizable systems.}
The proposed method is developed starting from a self-triggered implementation of the Lebesgue sampling rule. We recall that the Lebesgue sampling consists in updating the control law every time the triggering condition $\|x_{k}-x(t)\|\le \delta, \delta >0$ is violated,~\cite{ASTROM-CDC-2002}. 
Since, self-triggered sampling consists in predicting event occurrences, its design requires an upper-bound of the evolution of $\|x_{k}-x(t)\|$, which is given in the next result. 
\begin{lemma}\label{lem:upper-bound}
Let $M_{1}$ and $M_{2}$ be two positive constants such that the trajectories of~\eqref{eqn:closed-loop CT} satisfy
$
\|\xi(t)\|\le M_{1} \|\xi_{k}\|+M_{2}\bar w\,,
$

for all $t\ge t_k$. Then, the function $g(t):=x_{k}-x(t)$ is upper-bounded with

\begin{equation}\label{eqn:upper-bound}
\|g(t)\| \le (M_{1}\|x_{k}\|+M_{2}\bar w)(e^{L(t-t_{k})}-1)\,,
\end{equation}
where $L:=\max_{\eta\in \mathcal{D}_{\eta}} L_{f,u}L_{k,x}$,  for all $t\in[t_{k},t_{k+1})$.
\end{lemma}\\
A self-triggered sampler is devised by predicting the next time in which the function $\|g(t)\|$ hits the triggering threshold $\delta$, as done in~\cite{TIBERI-AUTO-2013} for linear systems.
This is equivalent to define $\Gamma(x,x_k)=\|g(t)\|-\delta$ and to predict the time instant $t_{k+1}$ for which $\Gamma(x,x_k)=0$ as per~\eqref{eqn:STC def}.
Such a prediction is performed by exploiting the bound~\eqref{eqn:upper-bound}, as stated in the the following result. 
\begin{proposition}\label{prop:Lebesgue sts}
Consider the same notation as in Lemma~\ref{lem:upper-bound} and let $\delta$ any arbitrary positive constant such that the trajectories of the perturbed system $\dot \psi=f(\eta,\psi,\kappa(\psi),w)+ Lg(t)$ where $\psi \in \mathcal D_{\xi}$ and $\| g(t)\|\le \delta$ are contained into the region of attraction $\mathcal{R}_a \subseteq \mathcal {D}_\xi$.
Then, the self-triggered sampler
\begin{equation}\label{eqn:Lebesgue sts}
t_{k+1}=t_{k}+\frac{1}{L}\ln \bigg(1+\frac{\delta  }{M_{1}\|x_k\|+M_{2}\bar w} \bigg)\,,
\end{equation}
ensures UUB of the sampled-data system~\eqref{eqn:sampled-data}. Moreover, there exists a positive constant $h_{\min}$ such that $t_{k+1}-t_{k}>h_{\min}$ for all $k$.
\end{proposition}\\
The self-triggered sampler~\eqref{eqn:Lebesgue sts} applies to every exponentially stable system of the form~\eqref{eqn:closed-loop CT}, and it is robust with respect to parameter uncertainty and external disturbances.
As~\eqref{eqn:Lebesgue sts} suggests, the inter-sampling intervals increase as $\delta$ does, but, on the other hand, the size of the ultimate bound also increases since the perturbation due to the sampling exhibits larger amplitudes. This means that $\delta$ can be intended as a tuning parameter that encodes the trade-off between inter-sampling intervals and ultimate-bound size. While the self-triggered sampler~\eqref{eqn:Lebesgue sts} presents only a single tuning parameter,  the following result provides more flexibility, since it allows the tuning of few more parameters.   
\begin{theorem}\label{thm:main}
Consider the same notation and assumptions of Proposition~\ref{prop:Lebesgue sts}, and let $\nu_0,\nu_1,\nu_3$ and $\nu_2$ arbitrary positive constants such that 
\begin{equation}\label{eqn:tuning rule sim}
\max_{\substack{r\ge 0}}  (M_1r+M_2\bar d)L\left[\left(1+\frac{\nu_0}{\nu_2r+\nu_3} \right)^{\frac{L}{\nu_1}}-1\right]\le \delta\,.
\end{equation}
Then, the self-triggered sampler
\begin{equation}\label{eqn:universal sts}
t_{k+1}=t_{k}+\frac{1}{\nu_1}\ln \bigg(1+\frac{\nu_0 }{\nu_2\|x_k\|+\nu_3} \bigg)\,,
\end{equation}
ensures UUB of the sampled-data system~\eqref{eqn:sampled-data}. Moreover, there exists a positive constant $h_{\min}$ such that $t_{k+1}-t_{k}>h_{\min}$ for all $k$.
\end{theorem}\\
%
In addition to the size of the ultimate bound, the self-triggered sampler~\eqref{eqn:universal sts} also allows to regulate the minimum and the maximum inter-sampling intervals. For instance, let us define $h_k:=t_{k+1}-t_k$ and 
\begin{align}
h_{\min}&:=\frac{1}{\nu_1}\ln \bigg(1+\frac{\nu_0 }{\nu_2(M_1\|x_0\|+M_2)+\nu_3} \bigg) \,, \label{eqn:h_min}\\
h_{\textrm{mid}}&:=\frac{1}{\nu_1}\ln \bigg(1+\frac{\nu_0 }{\nu_2b+\nu_3} \bigg) \,,\label{eqn:h_mid}\\
h_{\max}&:=\frac{1}{\nu_1}\ln \bigg(1+\frac{\nu_0 }{\nu_3} \bigg) \,.\label{eqn:h_max}
\end{align}
For all $k>0$ it holds $h_k\in [h_{\min},h_{\max}]$, while  for a sufficiently large $k'$ it holds $h_k\in [h_{\textrm{mid}},h_{\max}]$ for all $k>k'$. This means that for $k>k'$ the closed-loop system can be viewed as a periodically sampled system with period $h^*$ and jitter $\tilde h_k$ such that $h^*\pm \tilde h_k \in [h_{\textrm{mid}},h_{\max}]$. If the system~\eqref{eqn:sampled-data} with $w=0$ is exponentially stable for any time-varying $h_k \in [h_{\textrm{mid}},h_{\max}]$, then the self-triggered sampler~\eqref{eqn:universal sts} ensures exponential stability and not only UUB. This fact suggests a tuning method: for instance, it is enough to compute first the values of $h_{\textrm{mid}}$ and $h_{\max}$ to ensure exponential stability of~\eqref{eqn:sampled-data}, and then to tune the coefficients $\nu_i$'s according to~\eqref{eqn:h_min}--\eqref{eqn:h_max}.
As an example, one can consider the method in~\cite{SEURET-CDC-2009}, \cite{BRIAT-SCL-2012} to compute $h_{\textrm{mid}}$ and $h_{\max}$ for the linear case or~\cite{HU-AUTO-2000} for the nonlinear case. 
Furthermore, the inter-sampling intervals asymptotically converge to a constant sampling period $h^*=h_{\max}$. 
%
%
%
%
%
%
%
%
%
%
%
\subsection{Local analysis: asymptotically stabilizable systems.}
In the previous section we presented a self-triggered formula which applies to \emph{every} exponentially stabilizable system. In this section the results are extended to asymptotically stabilizable systems. Unfortunately, the nice property of~\eqref{eqn:universal sts} that allows its utilization with every exponentially stabilizable systems has no counterpart for asymptotically stabilizable systems. In-fact, a key ingredient to obtain the self-triggered formula~\eqref{eqn:universal sts} relies on the bound $\|\xi(t)\|\le M_{1}\|\Blue{\xi_{k}}\|+M_{2}\bar w$ which applies to \emph{every} exponentially stable systems of the form~\eqref{eqn:closed-loop CT}. For asymptotically stabilizable nonlinear system, such a bound is replaced by $\|\xi(t)\| \le \beta(\|\xi_k\|,0)+\gamma(\bar w)$, where $\beta$ is a class-$\K\L$ function and $\gamma$ is an appropriate class-$\K$ function which depends on the system under exam.
 
\begin{theorem}\label{thm:asympt}
Assume that the system~\eqref{eqn:closed-loop CT} is ISS with respect to $w$.
Let $\delta$ any arbitrary positive constant such that the trajectories of the perturbed system $\dot \psi=f(\eta,\psi,\kappa(\psi),w)+ Lg(t)$ where $\psi \in \mathcal D_{\xi}$ and  $\| g(t)\|\le \delta$ satisfy $\|\psi(t)\| \le \beta(\|\psi_k\|,t)+\gamma_1(\bar w)+\gamma_2(\delta)$ for some class-$\mathcal{KL}$ function $\beta$ and some class-$\mathcal K$ functions $\gamma_1, \gamma_2$, and assume that they are contained into the region of attraction $\mathcal{R}_a \subseteq \mathcal {D}_\xi$.  Finally, let $\nu_0,\nu_1$ and $\nu_3$ be arbitrary positive constants and let $\nu_2:\R_0^+\to\R_0^+$ be any function such that 
\begin{equation}\label{eqn:tuning rule sim asympt}
\max_{\substack{r\ge 0}}  (\beta(r,0)+\gamma_1(\bar w))L\left[\left(1+\frac{\nu_0}{\nu_2(r)+\nu_3} \right)^{\frac{L}{\nu_1}}-1\right]\le \delta\,.
\end{equation}
Then, the self-triggered sampler
\begin{equation}\label{eqn:sts nonlinear}
t_{k+1}=t_{k}+\frac{1}{\nu_1}\ln \bigg(1+\frac{\nu_0 }{\nu_2(\|x_k\|)+\nu_3} \bigg)\,,
\end{equation}
ensures UUB of the sampled-data system~\eqref{eqn:sampled-data}. Moreover, there exists a positive constant $h_{\min}$ such that $t_{k+1}-t_{k}>h_{\min}$ for all $k$.
\end{theorem}\\

The design of~\eqref{eqn:sts nonlinear} does not require exact knowledge of the function $\beta$, but only its behavior with respect to its first argument. Since it is well known that for asymptotically stable systems there exists a Lypaunov function $V(\xi)$ such that $\alpha_1(\|\xi\|)\le V(\xi)\le\alpha_2(\|\xi\|)$ for some class-$\K$ functions $\alpha_1$ and $\alpha_2$, it is enough to set $\nu_2(r)=\alpha_1^{-1}(V(r))$.
A notable case is when the closed-loop system admits a quadratic Lyapunov function, for which it holds $\alpha_1^{-1}(V(r))=cr$ for some positive $c$. In this case, the self-triggered sampler~\eqref{eqn:sts nonlinear} reduces to~\eqref{eqn:universal sts}. 
\begin{corollary}\label{cor:asympt}
Consider the same assumptions as in Theorem ~\eqref{thm:asympt} and assume that the closed-loop system with continuous control~\eqref{eqn:closed-loop CT} admits a quadratic Lyapunov function. Then, the self-triggered sampler~\eqref{eqn:universal sts} ensures UUB of the sampled-data system~\eqref{eqn:sampled-data} and there exists a positive constant $h_{\min}$ such that $t_{k+1}-t_{k}>h_{\min}$ for all $k$.
\end{corollary}\\
The tuning rules presented in the previous section also applies, \emph{mutatis mutandis}, to the self-triggered sampler~\eqref{eqn:sts nonlinear} as long as the function $\nu_2\|x_k\|$ is replaced with $\nu_2(\|x_k\|)$ and the bound $\|\xi(t)\|\le M_1\|\xi_k\|+M_2\bar w$ is replaced with $\|\xi(t)\|\le \beta(\|\xi_k\|,0)+\gamma_1(\bar w)$.
%
%
%
%
%
\subsection{Global analysis.}
The presented self-triggered sampler has been developed by modeling the deviation of the piecewise control from the continuous control as an external perturbation on the system state $g(t)$. The duty of the sampling strategy is to keep such a perturbation bounded. Then, by exploiting some \emph{bounded input - bounded state} property of the closed-loop system with continuous control, UUB of the sampled-data system is proved. The perturbation due to the proposed sampling rule is given by the left-hand side of~\eqref{eqn:tuning rule sim} or~\eqref{eqn:tuning rule sim asympt}. At this point, one may draw the conclusion that whenever the bounded input - bounded state property of the system~\eqref{eqn:closed-loop CT} holds globally, then the self-triggered samplers~\eqref{eqn:universal sts}--\eqref{eqn:sts nonlinear} applies for any setting of the parameters $\nu_i's$. Unfortunately, this is only partly true. By taking a closer look at~\eqref{eqn:tuning rule sim}, we observe that its right hand side is bounded if, and only if the Lipschitz constant $L$ holds globally. In-fact, there are many cases (e.g. polynomial systems) in which the system is locally, but not globally Lipschitz continuous. Hence, the method would only apply to bounded regions, and UUB would be only semi-global. Moreover, the selection of $\delta$ shall ensure boundedness of the trajectories into the region in which the Lipschitz constant $L$ is computed.

To tackle these issues, we recall that a local Lipschitz constant $L_f$ of a function $f$ depends on the domain in which it is computed. Let now $\mathcal F$ be the set of all the Lipschitz functions and let $Lip:2^{\R^n} \times \mathcal F \to \R_0^+$ be the operator that associates local Lipschitz constants $L_f$ to Lipschitz functions $f \in \mathcal F$ over subsets $\mathcal D\subseteq \R^n$. By observing that the proposed self-triggered sampler enforces the sets $\mathcal{B}_{d_k}$ to be invariant for all $t\in [t_k,t_{k+1})$, where $d_k=\beta(\|x_k\|,0)+\gamma_1(\bar w)+\gamma_2(\delta)$, it is not difficult to argue that there exists a function $\hat L:\R_0^+ \to \R_0^+$ such that $Lip(\mathcal{B}_{d_k}) \le \hat L(\|x_k\|)$ for all $t\in [t_k,t_{k+1})$ and for all $k$ as stated in the next result. 
\begin{lemma}\label{lem:1}
Let $\delta$ any arbitrary positive constant such that the trajectories of the perturbed system $\dot \psi=f(\eta,\psi,\kappa(\psi),w)+ Lg(t)$ with $\| g(t)\|\le \delta$ are UUB in any subset of $\R^{n_{\xi}}$ . 
Let $\mathcal F$ be the set of all the Lipschitz continuous functions and let $Lip:2^{\R^n} \times \mathcal F \to \R_0^+$ be the operator that associates local Lipschitz constants $L_f$ to Lipschitz continuous functions $f \in \mathcal F$ over subsets $\mathcal D\subseteq \R^n$. Then, by using the self-triggered sampler~\eqref{eqn:sts nonlinear}, there exists a function $\hat L:\R_0^+ \to \R_0^+$ such that $Lip(\mathcal{B}_{d_k}) \le \hat L(\|x_k\|)$ for all $t\in [t_k,t_{k+1})$ and for all $k$.
\end{lemma}\\
The next result represents a generalization of the self-triggered samplers presented in this note. 
\begin{theorem}\label{thm:global}
Consider the same notation as in Theorem~\ref{thm:asympt} and in Lemma~\ref{lem:1} and assume that the closed-loop system with continuous control~\eqref{eqn:closed-loop CT} is GUUB. Let $\nu_i:\R_0^+\to\R_0^+, i=1,\dots,3$ be any functions such that
\begin{equation}\label{eqn:tuning rule sim global}
\max_{\substack{r\ge 0}}  (\beta(r,0)+\gamma_1(\bar w))\hat L(r)\left[\left(1+\frac{\nu_0(r)}{\nu_2(r)+\nu_3(r)} \right)^{\frac{\hat L(r)}{\nu_1(r)}}-1\right]\le \delta\,.
\end{equation}
is bounded. Then, the self-triggered sampler
\begin{equation}\label{eqn:sts global}
t_{k+1}=t_{k}+\frac{1}{\nu_1(\|x_k\|)}\ln \bigg(1+\frac{\nu_0(\|x_k\|) }{\nu_2(\|x_k\|)+\nu_3(\|x_k\|)} \bigg)\,,
\end{equation}
ensures GUUB of the sampled-data closed-loop system. Moreover, there exists a positive constant $h_{\min}$ such that $t_{k+1}-t_{k}>h_{\min}$ for all $k$.
\end{theorem}\\
Note that the functions $\nu_i$ are not required to have any particular form. Infact, their duty is simply to bound the term in the left-hand side of~\eqref{eqn:tuning rule sim global}. In the next section we revise the example of Section~\ref{sec:motivating example} with our method.
\section{Motivating example revisited}\label{sec:simulations}
\begin{figure}[tp]
\centering
   \includegraphics[scale=0.63]{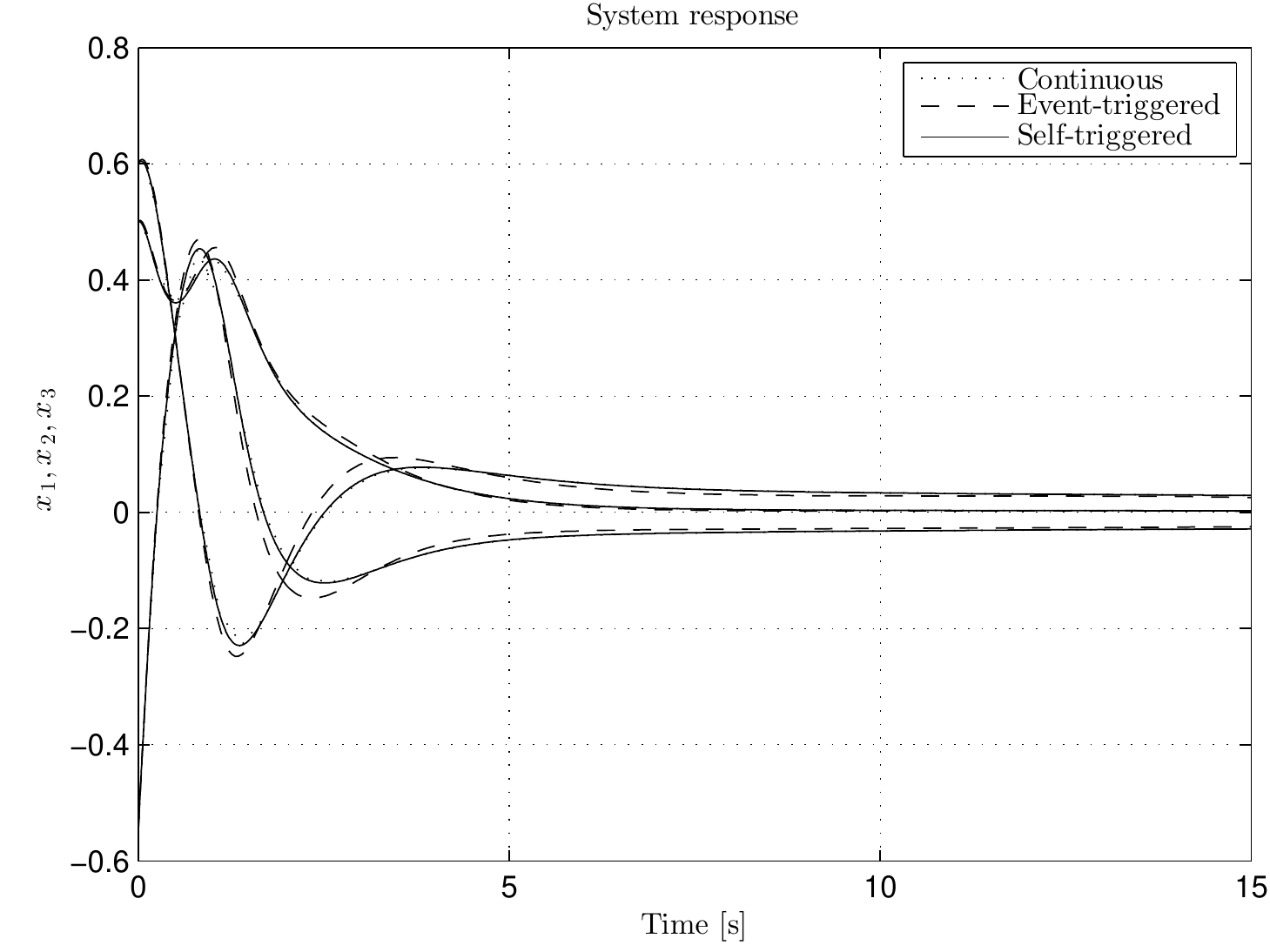}
   \caption{System response with $w=0$.}
   \label{fig:EvoUBA} 
\end{figure}
\begin{figure}[tp]
\centering
   \includegraphics[scale=0.63]{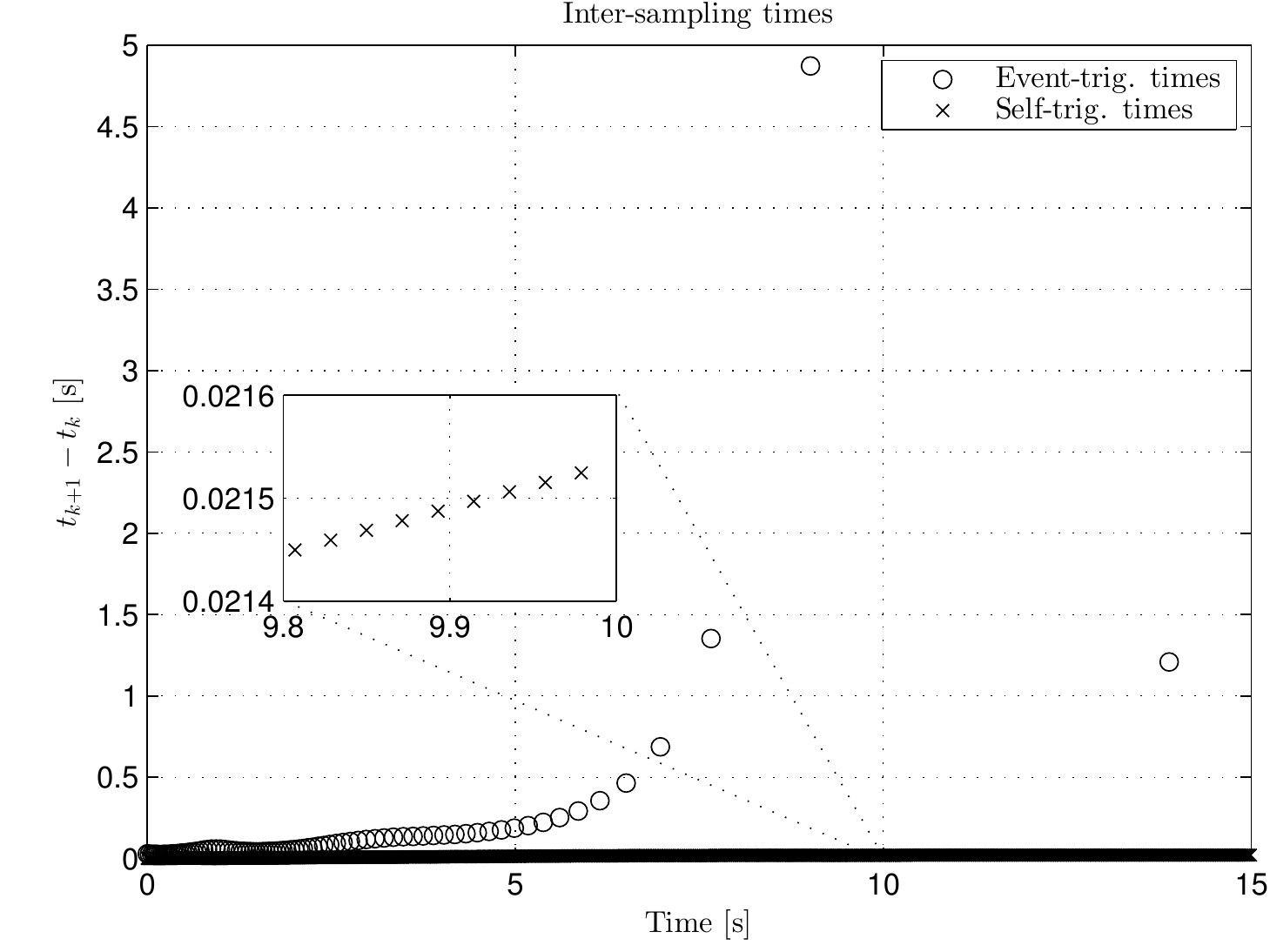}
   \caption{Inter-sampling times with $w=0$.}
   \label{fig:TrigUBA}
\end{figure}
%
%
%
%
%
\begin{figure}[tp]
\centering
   \includegraphics[scale=0.63]{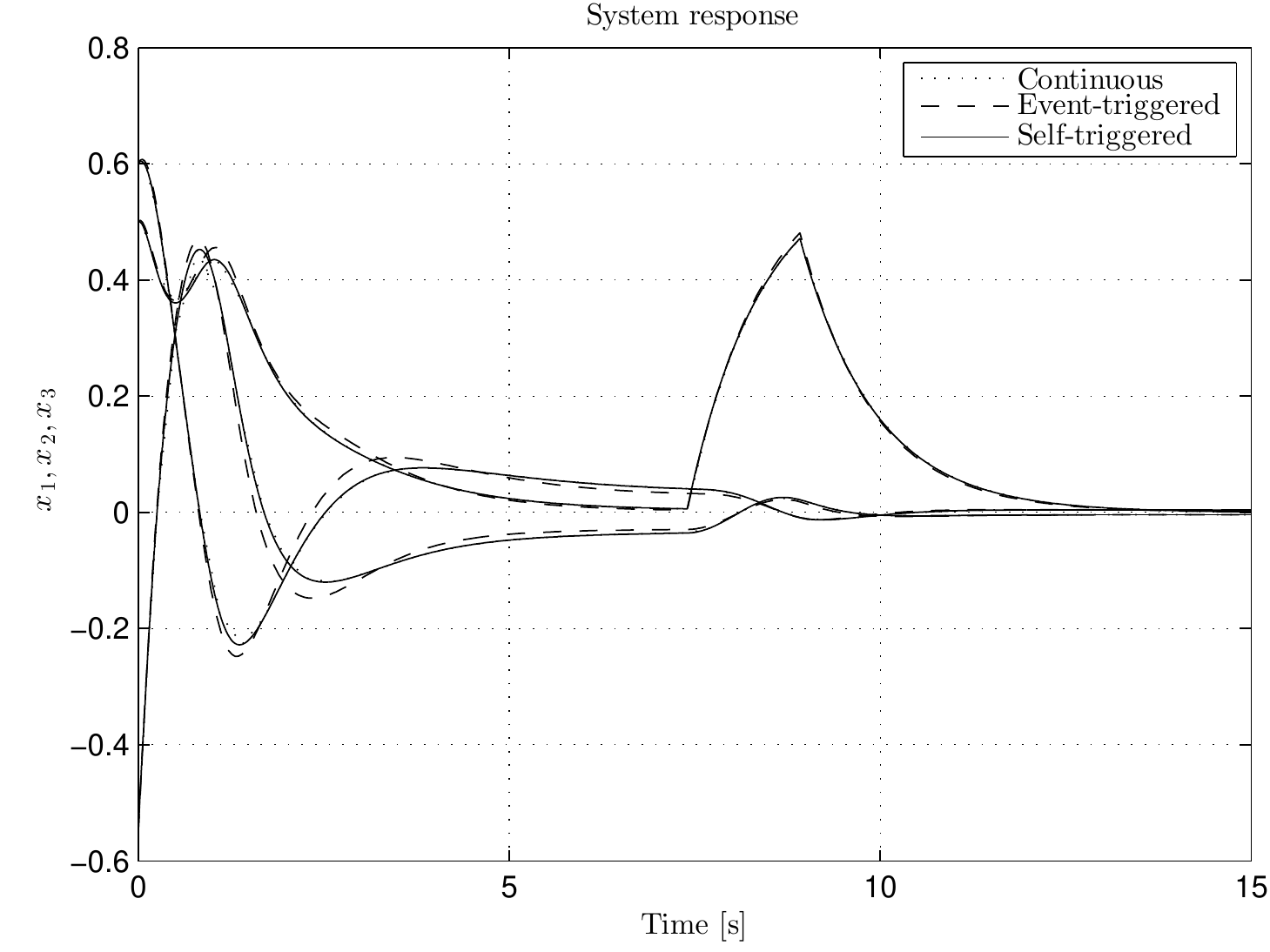}
   \caption{System response with $w=0.6$.}
    \label{fig:EvoUBAW}
\end{figure}
\begin{figure}[tp]
\centering
   \includegraphics[scale=0.63]{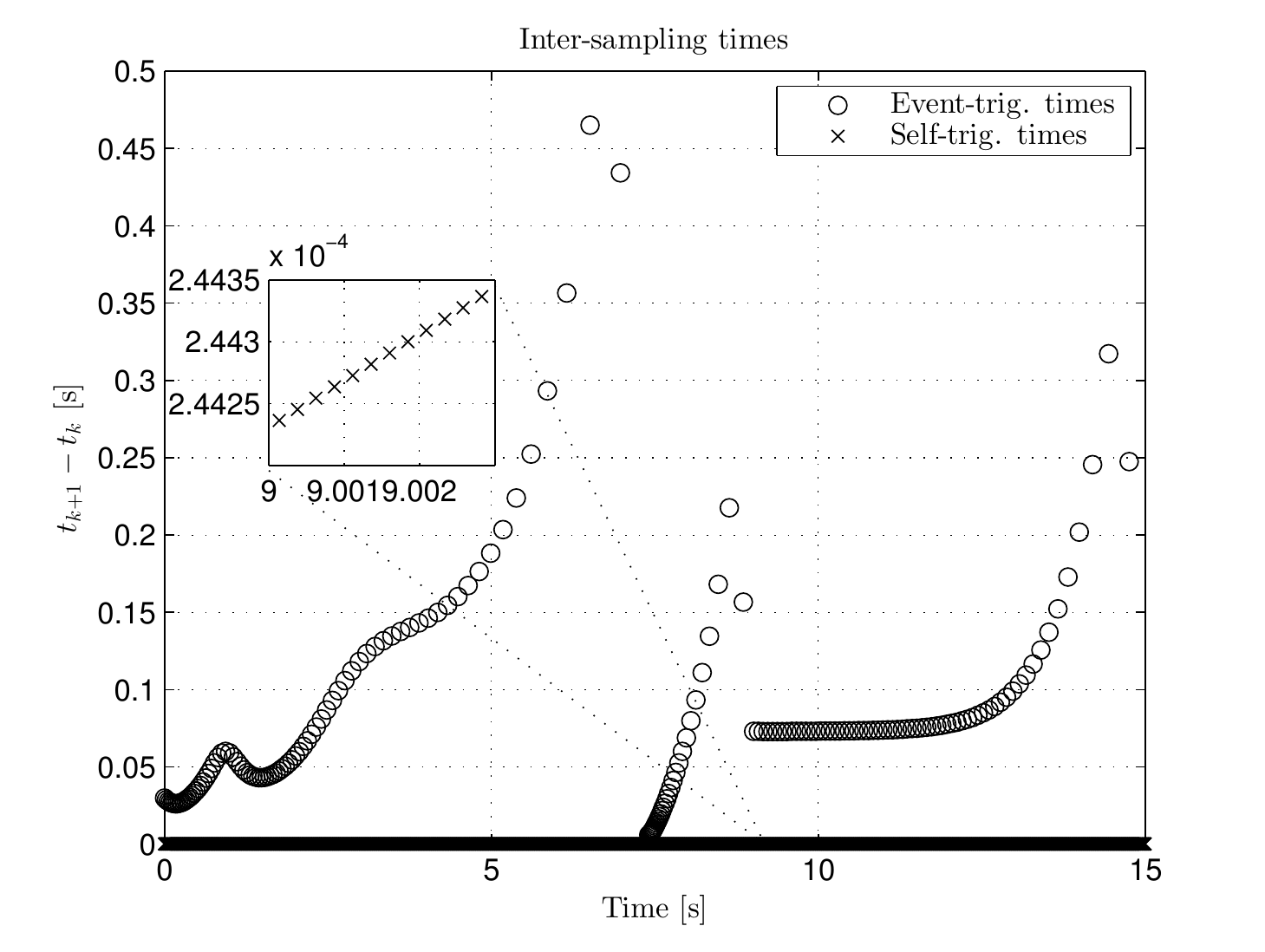}
   \caption{Inter-sampling times with $w=0.6$.}
   \label{fig:TrigUBAW}
\end{figure}
%

In this section we apply our method to the rigid-body control example described in Section~\ref{sec:motivating example}. We also compare our approach with existing methods by means of conservativeness of the inter-sampling intervals and closed-loop performance. The conservativeness of the inter-sampling times is measured in terms of average sampling time, while the closed-loop performance is evaluated through the average value $J_\textrm{avg}$ of a quadratic performance $J$ given by
\begin{equation*}
J=\int_0^{T}\|x(\sigma)\|^2+\|u(\sigma)\|^2\,d\sigma\,.
\end{equation*}
We consider 25 initial conditions equally spaced on a ball of radius one and a simulation time $T= 15$ s.  
Using SOSTOOLS~\cite{SOSTOOLS} we get a local parameter-dependent Lyapunov function that satisfies $3.2107\cdot 10^{-12}\|\xi\|^2 \le V(\eta,\xi) \le 3.6862 \cdot 10^{-12}\|\xi\|^2, \dot V \le -2.1273\cdot 10^{-13}\|\xi\|^4, \|\partial V/\partial \xi\| \le 2.34\cdot 10^{-13}\|\xi\|$ for all $\eta \in [1,8.2]$ on a ball of radius 5. We further get $L=61.1945$ over such a region. By setting $\delta = 2.8$ we get $b=3.9633$ as ultimate bound, and by considering $\mathcal B_1$ as initial condition set, we get $\|\psi\|_{\mathcal L_{\infty}}\le 4.7982$. 
Within this setting, we tune the proposed self-triggered sampler with $\nu_0(r) = 0.15\delta, \nu_1(r)=L, \nu_2(r)=10r$ and $\nu_3=10^{-6}$, for which we get $h_{\min}=0.142$ ms, $h_{\textrm{mid}}=0.180$ ms and $h_{\max}=211.6$ ms. 

The simulation results are reported in Table~\ref{tab:tab1}, while the system response for a particular initial condition is depicted in Figures~\ref{fig:EvoUBA}--\ref{fig:TrigUBA}. The system response with the self-triggered sampler~\cite{ANTA-TAC-2012} provides the larger inter-sampling times, but also the worst performance. Other methods exhibits comparable performance compared to the continuous-time case, meaning that the parameters uncertainty is well handled. However, compared to the other methods providing the same performance, our approach provides the largest average inter-sampling intervals.

Next, we evaluate the robustness with respect to external disturbances. For this purpose, we modify the dynamics of~\eqref{eqn:rigid-body} by setting $\dot \xi_2 = u_2+w$, where $w=0.6$ is an external disturbance acting on the system for $t \in [7.4,8.92]$ s. Within this setting, we set $\nu_0(r) = 0.05\delta, \nu_1(r)=L, \nu_2(r)=3.3r$ and $\nu_3=7.86$, for which we get $h_{\min}=0.009$ ms, $h_{\textrm{mid}}=0.048$ ms and $h_{\max}=0.288$ ms. 
As shown in Figure~\ref{fig:EvoUBAW} and as reported in Table~\ref{tab:tab2}, the response of the continuous-time, the event and proposed self-triggered implementation of the controller is fairly similar, whereas the method in~\cite{ANTA-TAC-2012} provides a deterioration of the performance index $J_{avg}$. A good trade-off between performance and average inter-sampling is instead provided by~\cite{ANTA-TAC-2010}. Nevertheless, there are no rigorous proofs of its robustness. Just for the sake of comparison, we tuned our self-triggered sampler by assuming $w=0$ even when in reality it holds $w = 0.6$.  With this setting, there are no proof of robustness of our method as well. Nevertheless, we experienced $J_{avg}=3.9413$ with an average sampling interval of $3.4$ ms, which outperforms~\cite{ANTA-TAC-2010}

Notice that by comparing our method in the perturbed and unperturbed case, the former exhibits a worsening of the average sampling period.  
This is due because in the tuning of the proposed self-triggered sampler when $w \neq 0$, a worst case disturbance acting for all the time has been considered. A method to reduce such conservativeness resorts to the utilization of disturbance observers as described in~\cite{TIBERI-NAHS-2013}. Nevertheless, differently from~\cite{TIBERI-NAHS-2013}, here we are dealing with an uncertain sampled-data nonlinear system, and to the best of our knowledge there are no result yet related to observer for this class of systems. Finally, for $w=0.6$, our self-triggered sampler provides an ultimate bound $b=3.9624$, whereas the other methods provide an ultimate bound $b=2.5682$. The larger ultimate bound in our case is due to the perturbation due to the sampling than sums up to $w$, while in the other cases, the ultimate bound only depends on the disturbance upper-bound $\bar w$. 

%
\begin{table*}[tp]
\centering
\resizebox{\textwidth}{!}{
\begin{tabular}{ccccccc}
\hline
 &Continuous & Event-trig. & Self-trig.~\cite{ANTA-TAC-2012} &Self-trig.~\cite{ANTA-TAC-2010}& Self-trig.~\cite{TIBERI-NAHS-2013}& Proposed Self-trig. \\
\hline
$J_\textrm{avg}$ &3.3921&3.6959&4.3298&3.4107&3.400&3.4101\\
Avg. time [ms]&-&155.7&182.5&2.500&0.7803&10.8\\
\hline
\end{tabular}
}
\caption{Average sampling intervals and $J_{\mathrm{avg}}$ for all the considered cases when $w=0$.}
\label{tab:tab1}
\end{table*}
\begin{table*}[tp]
\centering
\resizebox{\textwidth}{!}{
\begin{tabular}{ccccccc}
\hline
 &Continuous & Event-trig. & Self-trig.~\cite{ANTA-TAC-2012} &Self-trig.~\cite{ANTA-TAC-2010}& Self-trig.~\cite{TIBERI-NAHS-2013}& Proposed Self-trig. \\
\hline
$J_\textrm{avg}$ &3.9212&4.2582&5.4024&3.9523&3.9411&3.9364\\
Avg. time [ms]&-&77.1&93.4&2.500&0.7803&0.2678\\
\hline
\end{tabular}
}
\caption{Average sampling intervals and $J_{\mathrm{avg}}$ for all the considered cases when $w=0.6$.}
\label{tab:tab2}
\end{table*}
\section{Conclusions}
In this note we addressed the problem of robustness with respect to model uncertainties of self-triggered sampling for nonlinear systems. 
We have shown that even if a continuous-time controller is robust, this is not sufficient to use an arbitrary self-triggered sampling scheme, but the employed sampling scheme shall be also robust.  
In case of perfect model knowledge, then the self-triggered sampler in~\cite{ANTA-TAC-2012} outperforms our method, whereas in case of model uncertainties, our method appears to be more robust. 

A notable characteristic of the self-triggered sampler~\eqref{eqn:universal sts} relies in its similarity to existing methods,~\cite{TIBERI-AUTO-2013},\cite{LEMMON-TAC-09},\cite{MAZO-AUTOMATICA-2010},\cite{TOLIC-MED-2012}. This means that the proposed self-triggered formula can be regarded as a  generalization of existing self-triggered sampler that in turn can be used to tune~\eqref{eqn:universal sts} by  matching the coefficients $\nu_i$. For example, in the case of Lebesgue sampling~\eqref{eqn:Lebesgue sts}, a tuning rule is given by is  $\nu_0= \delta, \nu_1=L,\nu_2=M_1$ and $\nu_3= M_{2}\bar w$. The tuning of the parameter $\nu_3$ can also be performed through the utilization of a disturbance observer, as motivated and explained in~\cite{TIBERI-NAHS-2013}. During the process of coefficient matching, eventual parameter uncertainties can be easily included. Furthermore, the computational complexity of the proposed method is very low, since the next sampling-instant can be entirely determined by evaluating a simple function and no numerical methods are required.

Finally, we wish to highlight that in the system architecture definition we assumed that the self-triggered sampler is implemented in the controller side. However, our method still applies even whenever the self-triggered sampler is implemented on the sensor side as long as the controller updates are performed correspondingly to the output measurement transmission times. Eventual time delays can be easily accommodated by following the same line as in~\cite{TIBERI-AUTO-2013}.

%
%
\appendix
\begin{proof}{Proof of Lemma~\ref{lem:upper-bound}.}
First of all, note that it holds $ \dot g(t)=-\dot x(t)$ and $g(t_{k})=0$ at each sampling instant. For $t\in [t_{k},t_{k+1})$, it holds
\begin{align}
g(t)=&\int_{t_{k}}^{t}-f(\eta(s),x(s),\kappa(x_{k}),d(s))\,ds \\
=&\int_{t_{k}}^{t}f(\eta(s),x(s),\kappa(x(s)),d(s))\nonumber \\
&-f(\eta(s),x(s),\kappa(x_{k}),d(s))\,ds \nonumber \\
&-\int_{t_{k}}^{t}f(\eta(s),x(s),\kappa(x(s)),d(s))\,,
\end{align}
By taking the norm at both sides, and by recalling that exponential stability of~\eqref{eqn:closed-loop CT} with $w=0$ ensure the existence of constants $M_1,M_2$ and $\lambda$ such that $\|\xi(t)\|\le M_1\|\xi_k\|e^{-\lambda(t-t_k)}+M_2\bar w$ for all $t\ge t_k$, it follows
\begin{align}
\|g(t)\| \le& \int_{t_{k}}^{t}L_{f,u}L_{\kappa,x}\|x_{k}-x(s)\|\,ds \nonumber \\
&+\left\| \int_{t_{k}}^{t}f(\eta(s),x(s),\kappa(x(s)),d(s))\,ds \right\|\, \\
\le& \int_{t_{k}}^{t}L_{f,u}L_{\kappa,x}\|g(s)\|\,ds+M_{1}\|x_{k}\|+M_{2}\bar w\,\nonumber \\
\le & \int_{t_{k}}^{t}L\|g(s)\|\,ds+M_{1}\|x_{k}\|+M_{2}\bar w\,
\end{align}
By applying the Gronwall-Bellman inequality, it follows~\eqref{eqn:upper-bound}. 
\end{proof}
\begin{proof}{Proof of Proposition~\ref{prop:Lebesgue sts}.}
Converse Theorems ensure the existence of a parameters-dependent Lyapunov function $V(\eta,x)$\footnote{In case of time-varying uncertainty, the Lyapunov function $V(\eta,x)$ shall be replaced with a parameters-independent Lyapunov function $V(x)$, see~\cite{AMATO-BOOK}.}  that satisfies

\begin{equation}
\begin{aligned}
c_1\|\xi\|^2\le V(\eta,\xi)
&\le c_2\|\xi\|^2\,, \cr
\frac{\partial V(\xi)}{\partial \xi}f(\eta,\xi,\kappa(\xi),d)&\le-c_3\|\xi\|^2+c_5\bar w\,, \cr
\bigg\|\frac{\partial V(\eta,\xi)}{\partial \xi}\bigg\|&\le c_4\|\xi\| \,,  \cr
\end{aligned}   \label{eqn:LyapunovCandidateNonlinear}  
\end{equation}

where $c_1,c_2,c_3,c_4$ and $c_5$ are positive constants.
The time derivative of $V$ along the trajectories of the sampled-data system~\eqref{eqn:sampled-data}, for $t \in (t_{k},t_{k+1})$, satisfy
\begin{align}
\dot V &\le -c_3\|x\|^2+c_5\bar w+c_4\|x\|L\|x_k-x\| \nonumber \\
&\le -c_3\|x\|^2+c_5\bar w+c_4\|x\|L\hat g(x_k,t)\,,
\end{align}
for all  $t \in (t_{k},t_{k+1})$, where $\hat g(x_k,t)=(M_{1}\|x_{k}\|+M_{2}\bar w)(e^{L(t-t_{k})}-1)$. By setting $\hat g(x_k,t)=\delta$ it follows~\eqref{eqn:Lebesgue sts}. Moreover, since the sampling rule~\eqref{eqn:Lebesgue sts} enforces  $\|g(t)\|\le \delta$ for $t\in(t_k,t_{k+1})$, it follows $\|x(t)\|\le M_1\|x_k\|+M_2\bar w+M_3\delta:=M(k)$. Hence, the Lyapunov derivative is further upper-bounded with 
$\dot V\le -c_3\|x\|^2+c_5\bar w+c_4M(k)L\delta \,.$
By observing that at each sampling instants $t=t_k$ it holds $\dot V \le -c_3\|x_k\|^2+c_5\bar w$, it follows that the trajectories are upper-bounded for all $t\in [t_k,t_{k+1})$. Finally, since the sampling rule enforces $M(k+1)\le M(k)$ for $\|x\| > b$, or $M(k) \le b$ for $\|x\| \le b$ the sampled-data system~\eqref{eqn:sampled-data} is UUB.

\begin{equation}
h_{\min}:=\frac{1}{L}\ln\left(1+\frac{\delta}{M_1 c+M_2\bar w}\right)\,.
\end{equation}

where $c = \max\{M(0),b\}$.

Next, we have to prove that the sampling rule~\eqref{eqn:Lebesgue sts} guarantees that $x(t)\in \mathcal{R}_a$ for all $t\ge t_0$. Since by assumption the trajectories $\psi$ are confined into the region of attraction $\mathcal{R}_a$ for all $t \ge t_0$, and since~\eqref{eqn:Lebesgue sts} enforces $\|g(t)\|\le \delta$ and since $\lim_{t\to t_k^+} x(t)= \lim_{t\to t_k^-} x(t)$ for all $k$, it follows that $x(t)\in \mathcal{R}_a$ for all $t\ge t_0$. 

\end{proof}
%
%
%
%
%
%
%
%
%
%
\begin{proof}{Proof of Theorem~\ref{thm:main}.}
The first part of the proof follows the same line as the proof of Proposition~\ref{prop:Lebesgue sts}. This means that the Lyapuonov function $V(x)$ satisfies
$\dot V \le -c_3\|x\|^2+c_5\bar w+c_4\|x\|L\hat g(x_k,t) \,,$ where $\hat g(x_k,t)=(M_{1}\|x_{k}\|+M_{2}\bar w)(e^{L(t-t_{k})}-1)$.
By using the sampling rule~\eqref{eqn:universal sts}, it follows that
\begin{align}
\dot V\le& -c_3\|x\|^2+c_5\bar w+c_4\|x\|L \nonumber \\ 
&\times (M_1\|x_k\|+M_2\bar w)\left[\left(1+\frac{\nu_0}{\nu_2\|x_k\|+\nu_3} \right)^{\frac{L}{\nu_1}}-1\right] \nonumber \\
\le& -c_3\|x\|^2+c_5\bar w+c_4\|x\|L\bar \delta \,,
\end{align}
for all $t\in (t_k,t_{k+1})$, where
\begin{equation}
\bar \delta :=\max_{r \in \R^+_0}(M_1r+M_2\bar w)\left[\left(1+\frac{\nu_0}{\nu_3r+\nu_2} \right)^{\frac{L}{\nu_1}}-1\right] \,.
\end{equation}
Since $\bar \delta$ is bounded, UUB follows for all $t\in (t_k,t_{k+1})$. Moreover, correnspondently to the sampling instants $t=t_k$ it holds $\dot V \le -c_3\|x\|^2+c_5\bar w$ and thus UUB holds in every interval of the form $[t_k,t_{k+1})$. Finally, continuity of the solution of the sampling-data system~\eqref{eqn:sampled-data} ensure UUB for all $[t_k,t_{k+1})$ and for all $k$.
The existence of a lower-bound of the inter-sampling intervals can be proved by using the same arguments as in the proof of  Proposition~\ref{prop:Lebesgue sts}.
\end{proof}
\begin{proof}Proof of Theorem~\ref{thm:asympt}.
By observing that in the case of asymptotic stability the bound $\|\xi(t)\|\le M_{1}\|\xi_{k}\|+M_{2}\bar w$ becomes $\|\xi(t)\|\le \beta(\|\xi_k\|,0)+\gamma_1(\bar w)$, the upper-bound~\eqref{eqn:upper-bound} becomes  
$\|g(t)\|\le (\beta(\|x_{k}\|,0)+\gamma_1(\bar w))(e^{L(t-t_{k})}-1)\,.$
By using such a bound, the proof follows analougsly to the proof of Theorem~\ref{thm:main}, where the terms $M_1\|\xi_k\|$ and $M_2\bar w$ are replaced with $\beta(\|\xi_k\|,0)+\gamma_1(\bar w)$, and where the comparison functions in~\eqref{eqn:LyapunovCandidateNonlinear} are replaced with appropriate class-$\K$ functions~\cite{KHALIL}.
\end{proof}

\bibliographystyle{IEEEtran}
\bibliography{IEEEabrv,./MyBib}

\end{document}